\documentclass[11pt]{article}
\usepackage{amsfonts}
\usepackage{mathrsfs}
\usepackage{amsfonts,amssymb}
\usepackage{amsmath,amscd}
\usepackage[all]{xy}
\usepackage{wasysym, stmaryrd}
\usepackage[dvips]{color}
\usepackage[dvipdfm,
            colorlinks=ture,
            citecolor=blue,
            linkcolor=black]{hyperref}
\begin{document}
\sloppy
\newcommand{\dickebox}{{\vrule height5pt width5pt depth0pt}}
\newtheorem{Def}{Definition}[section]
\newtheorem{Bsp}{Example}[section]
\newtheorem{Prop}[Def]{Proposition}
\newtheorem{Theo}[Def]{Theorem}
\newtheorem{Lem}[Def]{Lemma}
\newtheorem{Koro}[Def]{Corollary}
\newtheorem{Claim}[Def]{Claim}
\newcommand{\heiti}{\bf\CJKfamily{hei}}
\newcommand{\kaiti}{\bf\CJKfamily{kai}}

\newcommand{\lra}{\longrightarrow}
\newcommand{\ra}{\rightarrow}
\newcommand{\F}{\mathcal {F}}
\newcommand{\Hom}{{\rm Hom}}
\newcommand{\End}{{\rm End}}
\newcommand{\Ext}{{\rm Ext}}
\newcommand{\Tor}{{\rm Tor}}
\newcommand{\pd}{{\rm proj.dim}}
\newcommand{\inj}{{\rm inj}}
\newcommand{\lgd}{{l.{\rm gl.dim}}}
\newcommand{\gld}{{\rm gl.dim}}
\newcommand{\fd}{{\rm fin.dim}}
\newcommand{\Fd}{{\rm Fin.dim}}
\newcommand{\lfd}{l.{\rm Fin.dim}}
\newcommand{\rfd}{r.{\rm Fin.dim}}
\newcommand{\Mod}{{\rm Mod}}
\newcommand{\Proj}{{\rm Proj}}
\newcommand{\modcat}[1]{#1\mbox{{\rm -mod}}}
\newcommand{\pmodcat}[1]{#1\mbox{{\rm -proj}}}
\newcommand{\Pmodcat}[1]{#1\mbox{{\rm -Proj}}}
\newcommand{\injmodcat}[1]{#1\mbox{{\rm -inj}}}
\newcommand{\E}{{\rm E}^\Phi}
\newcommand{\X}{ \mathscr{X}_{\mathcal {F}}^{{\rm F},\Phi}}
\newcommand{\Y}{\rm \mathscr{Y}_{\mathcal {F}}^{{\rm F},\Phi}}
\newcommand{\A}{\mathcal {A}}
\newcommand{\C}{\rm \mathscr{C}}
\newcommand{\K}{\rm \mathscr{K}}
\newcommand{\D}{\rm \mathscr{D}}
\newcommand{\opp}{^{\rm op}}
\newcommand{\otimesL}{\otimes^{\rm\bf L}}
\newcommand{\otimesP}{\otimes^{\bullet}}
\newcommand{\rHom}{{\rm\bf R}{\rm Hom}}
\newcommand{\projdim}{\pd}
\newcommand{\stmodcat}[1]{#1\mbox{{\rm -{\underline{mod}}}}}
\newcommand{\Modcat}[1]{#1\mbox{{\rm -Mod}}}
\newcommand{\modcatr}[1]{\mbox{{\rm mod}#1}}
\newcommand{\Modcatr}[1]{\mbox{{\rm Mod}#1}}
\newcommand{\Pmodcatr}[1]{\mbox{{\rm Proj}#1}}
\newcommand{\procat}[1]{#1\mbox{{\rm -proj}}}
\newcommand{\Tr}{\rm Tr}
\newcommand{\add}{{\rm add}}
\newcommand{\Imf}{{\rm Im}}
\newcommand{\Ker}{{\rm Ker}}
\newcommand{\EA}{{\rm E^\Phi_\mathcal {A}}}
\newcommand{\pro}{{\rm pro}}
\newcommand{\Coker}{{\rm Coker}}
\newcommand{\id}{{\rm id}}
\renewcommand{\labelenumi}{\Alph{enumi}}
\newcommand{\M}{\mathcal {M}}
\newcommand{\Mf}{\rm \mathcal {M}^f}
\newcommand{\rad}{{\rm rad}}
\newcommand{\injdim}{{\rm inj.dim}}
{\Large \bf

\begin{center} Derived equivalences between subrings\end{center}}
\medskip
\centerline{\bf Yiping Chen }
\medskip \centerline{School of Mathematics and Statistics, Wuhan University, Wuhan 430072, China}
\medskip \centerline{School of Mathematical Sciences, Beijing Normal University, Beijing 100875, China}

\medskip
\centerline{E-mail: ypchen@whu.edu.cn}

\medskip
\abstract{
In this paper, we construct derived equivalences between two
subrings of relevant $\Phi$-Auslander-Yoneda rings from an arbitrary
short exact sequence in an abelian category. As a consequence, any
short exact sequence in an abelian category gives rise to a derived
equivalence between two subrings of endomorphism rings. These
results generalize some methods on constructing derived
equivalences. }



\section{Introduction}

Derived equivalences play an important role in many branches of
mathematics. For example, the representation theory of algebras,
combinatorics and algebraic geometry (see \cite{AAB, BGG, BK, LM,
LS}). One of fundamental problems in this field is how to construct
derived equivalences between rings. Recently, Hu and Xi have
exhibited derived endomorphism rings induced by $\mathcal{D}$-split
sequences \cite{HX1}. Let $\mathcal{C}$ be an additive category and
$\mathcal {D}$ be a full subcategory of $\mathcal{C}$. A sequence
$$X\stackrel{f}\ra M\stackrel{g}\ra Y\qquad(*)$$ is called a {\em $\mathcal {D}$-split sequence} if it
satisfies that $(1)$ $f$ is a kernel of $g$, and $g$ is a cokernel
of $f$; $(2)$ $\Hom_{\mathcal{C}}(f, M)$ and $\Hom_{\mathcal{C}}(M,
g)$ are onto; $(3)$ $M\in\mathcal {D}$ (see \cite[Definition
3.1]{HX1}). Given a $\mathcal {D}$-split sequence $(*)$, Hu and Xi
proved that the two rings $\End_{\mathcal{C}}(X\oplus M)$ and
$\End_{\mathcal{C}}(M\oplus Y)$ are derived equivalent \cite[Theorem
3.5]{HX1}. This result establishes a beautiful connection between
derived equivalences and $\mathcal {D}$-split sequences.

Let $\mathcal {A}$ be an abelian category. Our aim in this paper is
to extent this result and construct derived equivalences induced by
a short exact sequence in $\mathcal {A}$. Our ideal in this
direction is to find out an additive subcategory $\mathcal {B}$ of
$\mathcal{A}$ satisfying that a short exact sequence in
$\mathcal{A}$ is a $\mathcal {D}$-split sequence in $\mathcal{B}$.
In order to describe the main result precisely, we fix some notation
first. Let $\Phi$ be an admissible subset of natural numbers. We
denote the derived category of $\mathcal {A}$ by $D^b(\mathcal{A})$
with suspension functor $[1]$. Let $X$ be an object in
$D^b(\mathcal{A})$. The $\Phi$-Auslander-Yoneda algebra is
$\E_{\mathcal{A}}(X):=\oplus_{i\in\Phi}\Hom_{D^b(\mathcal{A})}(X,
X[i])$ with multiplication in a natural way. The matrix rings
$\begin{pmatrix}
  \widehat{\EA(X)}&\widehat{\EA(X,M)}\\
  \EA(M,X)&\EA(M)
\end{pmatrix}$ and $\begin{pmatrix}
\widehat{\EA(Y)}&\EA(Y,M)\\\widehat{\EA(M,Y)}&\EA(M)
\end{pmatrix}$ are some subrings of $\EA(X\oplus M)$ and $\EA(Y\oplus
M)$, respectively. For the concise definitions of these rings, we
refer the reader to section 3. Our result can be stated in the
following form:

\begin{Theo}\label{theorem} Let $\Phi$ be an admissible subset of
$\mathbb{N}$, and let $\mathcal {A}$ be an abelian category and $M$
an object in $\mathcal {A}$. Suppose that
$$0\ra X\stackrel{\alpha}\ra M_1\stackrel{\beta}\ra Y\ra 0$$
is an exact sequence in $\mathcal {A}$ satisfying that $M_1\in \add
(M)$ and $\Ext^i_\mathcal {A}(M, X)=0, \Ext^i_\mathcal {A}(Y, M)=0$
for $0\neq i\in \Phi$. Then the two rings
$$\begin{pmatrix}
  \widehat{\EA(X)}&\widehat{\EA(X,M)}\\
  \EA(M,X)&\EA(M)
\end{pmatrix}\mbox{ and }
\begin{pmatrix}
\widehat{\EA(Y)}&\EA(Y,M)\\\widehat{\EA(M,Y)}&\EA(M)
\end{pmatrix}$$
are derived equivalent.
\end{Theo}

Theorem \ref{theorem} partially extends the main result of Hu and Xi
in \cite[Theorem 3.5]{HX1}. As a corollary, any short exact sequence
in an abelian category implies a derived equivalence between two
subrings of relevant endomorphism rings.

\begin{Koro}Let $\mathcal {A}$ be an abelian category and $M$ an object in
$\mathcal{A}$. Suppose that
$$0\ra X\stackrel{\alpha}\ra M_1\stackrel{\beta}\ra Y\ra 0$$
is an exact sequence in $\mathcal{A}$ with $M_1\in\add(M)$.

Then the following two rings
$$\begin{pmatrix}
  \widehat{\End_{\A}(X)}&\widehat{\Hom_{\A}(X,M)}\\
  \Hom_{\A}(M,X)&\End_{\A}(M)
\end{pmatrix}\mbox{ and }  \begin{pmatrix}\widehat{\End_{\A}(Y)}
  &\Hom_{\A}(Y,M)\\\widehat{\Hom_{\A}(M,Y)}
  &\End_{\A}(M)
\end{pmatrix}$$
are derived equivalent, where
$$\widehat{\End_{\A}(Y)}:=\{t\in\End_{\A}(Y)\mid \mbox{there exists a morphism }t_1:{}M_1\ra M_1 \mbox{ such that }
\beta t=t_1\beta\mbox{ in } \A\}$$
$$\widehat{\End_{\A}(X)}:=\{t\in\End_{\A}(X)\mid \mbox{there exists a morphism }t_1:{}M_1\ra M_1 \mbox{ such that }
t\alpha=\alpha t_1 \mbox{ in } \A\}$$
$$\widehat{\Hom_{\A}(X,M)}:=\{t\in\Hom_{\A}(X,M)\mid \mbox{there exists a morphism }t_1:{}M_1\ra M \mbox{ such that }
t=\alpha t_1\mbox{ in } \A\}$$
$$\widehat{\Hom_{\A}(M,Y)}:=\{t\in\Hom_{\A}(M,Y)\mid \mbox{there exists a morphism }t_1:{}M\ra M_1 \mbox{ such that }
t=t_1\beta\mbox{ in } \A\}.$$
\end{Koro}

This paper is organized as follows: In section 2, we will recall
some notation and definitions. In section 3, we will give the proof
of Theorem \ref{theorem} which partially extends \cite[Theorem
3.5]{HX1}. As a byproduct of Theorem \ref{theorem}, we will point
out that any short exact sequence induces a derived equivalence
between two subrings of endomorphism rings. In particular, Corollary
\ref{D} generalizes Theorem 4.1 in \cite{D}. In section 4, we will
give an example to illustrate our results.
\section{Preliminaries}

In this section, we will recall some basic definitions and facts
which are needed in our proofs.

\subsection{Notation and conventions}

Throughout this paper, we assume that $\mathcal {A}$ is an abelian
category, and $\mathcal {T}$ is a triangulated category with
suspension functor $[1]$.

Let $\mathcal {C}$ be an additive category. For an object $M$ in
$\mathcal{C}$, we denote by $\add(M)$ the full subcategory of
$\mathcal{C}$ consisting of all direct summands of finite direct
sums of $M$. For two morphisms $f:{}X\ra Y$ and $g:{}Y\ra Z$ in
$\mathcal{C}$, we write $fg$ for their composition which is a
morphism from $X$ to $Z$. For two functors
$F:{}\mathcal{C}\ra\mathcal{D}$ and $G:{}\mathcal{D}\ra\mathcal{E}$,
we denote $G\circ F$ for their composition.

Let $\mathcal {A}$ be an abelian category. We denote by
$C(\mathcal{A})$ the category of complexes of objects in
$\mathcal{A}$. A complex $X^\bullet$ is called {\em bounded below}
if $X^i=0$ for all but finitely many $i<0$, {\em bounded above} if
$X^i=0$ for all but finitely many $i>0$, and {\em bounded} if
$X^\bullet$ is bounded above and bounded below. We denote by
$K(\mathcal{A})$ the homotopy category of $\mathcal{A}$, and denote
by $D(\mathcal{A})$ the derived category of $\mathcal{A}$. The full
subcategory of $K(\mathcal{A})$ and $D(\mathcal{A})$ consisting of
bounded complexes over $\mathcal{A}$ is denoted by
$K^b(\mathcal{A})$ and $D^b(\mathcal{A})$, respectively. We denote
by $C^-(\mathcal{A})$ the category of complexes of bounded above,
and by $K^-(\mathcal{A})$ the homotopy category of
$C^-(\mathcal{A})$. The full subcategory of $D(\mathcal{A})$
consisting of bounded above complexes is denoted by
$D^-(\mathcal{A})$. It is well known that the categories
$K^{-}(\mathcal{A}), K^b(\mathcal {A}), D^{-}(\mathcal{A})$ and
$D^b(\mathcal{A})$ are triangulated categories.

Suppose that all rings are associative with identity, and all
modules are unitary unless otherwise stated are left modules. Let
$A$ be a ring. We denote the category of all (respectively, finitely
generated) left $A$-modules by $\Modcat{A}$ (respectively,
$\modcat{A}$). The category of all projective left $A$-modules is
denoted by $\Pmodcat{A}$, and the category of all finitely generated
projective left $A$-modules is denoted by $\pmodcat{A}$.

\medskip
The following result is called the Morita theorem of derived
categories of rings.
\begin{Lem}$\cite{R1}$ Let $A$ and $B$ be two rings. The following
conditions are equivalent:

$(1)$ $K^-(\pmodcat{A})$ and $K^-(\pmodcat{B})$ are equivalent as
triangulated categories;

$(2)$ $D^b(\Modcat{A})$ and $D^b(\Modcat{B})$ are equivalent as
triangulated categories;

$(3)$ $K^b(\Pmodcat{A})$ and $K^b(\Pmodcat{B})$ are equivalent as
triangulated categories;

$(4)$ $K^b(\pmodcat{A})$ and $K^b(\pmodcat{B})$ are equivalent as
triangulated categories;

$(5)$ The two rings $B$ and $\End(T^\bullet)$ are isomorphic, where
$T^\bullet$ is a complex in $K^b(\pmodcat{A})$ satisfying:

\quad $(a)$ $T^\bullet$ is self-orthogonal, that is,
$\Hom_{K^b(\pmodcat{A})}(T^\bullet,T^\bullet[i])=0$ for all $i\neq
0$,

\quad $(b)$ $\add(T^\bullet)$ generates $K^b(\pmodcat{A})$ as a
triangulated category.
\end{Lem}

If the above conditions $(1)-(5)$ are satisfied, then the two rings
are called {\em derived equivalent}. The complex $T^\bullet$ in
$(5)$ is called a {\em tilting complex} over $A$.

Let $\mathbb{N}=\{0, 1, 2, \cdots\}$ be the set of natural numbers.
A subset $\Phi$ of $\mathbb{N}$ containing $0$ is called an {\em
admissible subset} of $\mathbb{N}$ if the following condition is
satisfied:

If $i, j$ and $k$ are in $\Phi$ such that $i+j+k\in\Phi$, then
$i+j\in\Phi$ if and only if $j+k\in\Phi$.

There are a lot of admissible subsets of $\mathbb{N}$. The sets
$\{0, 3, 4\}$ and $\{0, 1, 2, 3, 4\}$ are admissible subsets of
$\mathbb{N}$. If $\Phi$ is an admissible subset of $\mathbb{N}$,
then $m\Phi$ is an admissible subset of $\mathbb{N}$ for every
$m\in\Phi$. Suppose that $\Phi$ is a subset of $\mathbb{N}$
containing $0$, then the set $\Phi^m=\{x^m\mid x\in \Phi\}$ is an
admissible subset of $\mathbb{N}$ for all $m\geq 3$. For more
details, we refer the reader to \cite{HX2}.

Let $\Phi$ be an admissible subset of $\mathbb{N}$, and let
$\mathcal {T}$ be a triangulated category with suspension functor
$[1]$. Consider the {\em $\Phi$-orbit category} $\mathcal {T}^\Phi$
whose objects are the objects of $\mathcal {T}$. Suppose that $X, Y$
are two objects in $\mathcal {T}^\Phi$. The homomorphism set in
$\mathcal {T}^\Phi$ is defined as follows:
$$\Hom_{\mathcal {T}^\Phi}(X,Y):= \bigoplus_{i\in \Phi}\Hom_\mathcal {T}(X,Y[i])\in \Modcat{\mathbb{Z}}$$
and the composition is defined in an obvious way. Since $\Phi$ is an
admissible subset of $\mathbb{N}$, the $\Phi$-orbit category
$\mathcal {T}^\Phi$ is an additive category. Let $X, Y$ be objects
in $\mathcal {T}^\Phi$. The homomorphism set $\Hom_{\mathcal
{T}^\Phi}(X, X)$, denoted by $\E_{\mathcal {T}}(X)$, is a
$\mathbb{Z}$-algebra. It is called $\Phi$-Auslander-Yoneda algebra
of $X$ in \cite{HX2}. The homomorphism set $\Hom_{\mathcal
{T}^\Phi}(X, Y)$ is a $\E_\mathcal {T}(X)-\E_\mathcal
{T}(Y)$-bimodule. We denote it briefly by $\E_\mathcal {T}(X, Y)$.
In particular, if $\mathcal {T}=D^b(\mathcal {A})$ (respectively,
$D^b(\Modcat{A})$), then we abbreviate $\E_{D^b(\mathcal{A})}(X, Y)$
(respectively, $\E_{D^b(\Modcat{A})}(X, Y)$) to $\E_{\mathcal{A}}(X,
Y)$ (respectively, $ (\E_{A}(X, Y)$), where $X$ and $Y$ are objects
in $D^b(\mathcal{A})$ (respectively., $(D^b(\Modcat{A})$).

\medskip

The following Lemma which is taken from \cite{HX2}, can be verified
directly.

\begin{Lem}\cite[Lemma 3.5]{HX2}Let $\Phi$ be an admissible subset of $\mathbb{N}$. Let $A$ be a ring with identity and $V$ be an $A$-module.
Suppose that $V_1$ and $V_2$ are in $\add(V)$. Then the following
statements are true.

$(1)$ The $\E_A(V,V_1)$ is projective and finitely generated, and
there is an isomorphism
$$\mu: \E_A(V_1,V_2)\ra \Hom_{\E_A(V)}(\E_A(V,V_1),\E_A(V,V_2))$$
which sends $(f_i)\in \E_A(V_1,V_2)$ to the morphism $(a_i)\mapsto
(a_i)(f_i)$ for $(a_i)\in \E_A(V,V_1)$. Moreover, if $V_3\in\add(V)$
and $(g_i)\in \E_A(V_2,V_3)$, then
$\mu((f_i)(g_i))=\mu((f_i))\mu((g_i))$.

$(2)$ The functor $\E_A(V,-):{}\add(V)\ra \pmodcat{\E_A(V)}$ is
faithful.

$(3)$ If $V_1$ is projective or $V_2$ is injective, then the functor
$\E_A(V,-)$ induces an isomorphism of $\mathbb{Z}$-modules:
$$\E_A(V,-): \Hom_A(V_1,V_2)\ra \Hom_{\E_A(V)}(\E_A(V,V_1),\E_A(V,V_2)).$$
\end{Lem}

\section{Proof of Theorem \ref{theorem}}

In this section, we will construct derived equivalences from an
exact sequence in an abelian category $\mathcal{A}$. First, we will
prove Theorem \ref{theorem}, and then we will derive some useful
consequences from the main result.

Let $M$ be an object in $\mathcal {A}$. Suppose that $0\ra
X\stackrel{\alpha}\ra M_1\stackrel{\beta}\ra Y\ra 0\quad (**)$ is an
exact sequence in $\mathcal {A}$ with $M_1\in$ add($M$). By abuse of
notation, $(**)$ yields a triangle $X\stackrel{\alpha}\ra
M_1\stackrel{\beta}\ra Y \ra X[1]$ in $D^b(\mathcal {A})$. Let us
denote $X\oplus M\oplus Y$ by $V$. Consider the
$\Phi$-Auslander-Yoneda algebra
$$\EA(V):=\begin{pmatrix}
\EA(X)&\EA(X,M)&\EA(X,Y)\\
\EA(M,X)&\EA(M)&\EA(M,Y)\\
\EA(Y,X)&\EA(Y,M)&\EA(Y)
\end{pmatrix}.$$

For convenience, we denote $\EA(V)$ by $\Gamma$.

Select a subring $\Lambda$ of $\Gamma$
$$\Lambda:=\begin{pmatrix}
  \widehat{\EA(X)}&\widehat{\EA(X,M)}&\widehat{\EA(X,Y)}\\
\EA(M,X)&\EA(M)&\widehat{\EA(M,Y)}\\
\EA(Y,X)&\EA(Y,M)&\widehat{\EA(Y)}
\end{pmatrix}$$

where
$$\widehat{\EA(X)}:=\{(t_i)_{i\in \Phi}\in \EA(X)\mid t_i\alpha[i] \mbox{ factorizes through }\alpha \mbox{ in }
D^b(\mathcal {A}) \mbox{ for } i\in \Phi\},$$
$$\widehat{\EA(Y)}:=\{(t_i)_{i\in \Phi}\in\EA(Y)\mid\beta t_i \mbox{ factorizes through }\beta[i] \mbox{ in }
D^b(\mathcal {A}) \mbox{ for } i\in \Phi\},$$
$$\widehat{\EA(X,M)}:=\{(t_i)_{i\in\Phi}\in\EA(X,M)\mid t_i \mbox{ factorizes through }\alpha \mbox{ in }
D^b(\mathcal {A})\mbox{ for } i\in \Phi \},$$
$$\widehat{\EA(M,Y)}:=\{(t_i)_{i\in \Phi}\in \EA(M,Y)\mid t_i \mbox{ factorizes through }\beta[i] \mbox{ in }
D^b(\mathcal {A}) \mbox{ for } i\in \Phi\},$$
$$\widehat{\EA(X,Y)}:=\{(t_i)_{i\in \Phi}\in\EA(X,Y)\mid t_i \mbox{ factorizes through }\alpha
\mbox{ and }\beta[i] \mbox{ in }D^b(\mathcal {A}) \mbox{ for }
i\in\Phi\}.$$

One can check directly that the matrix rings $\begin{pmatrix}
  \widehat{\EA(X)}&\widehat{\EA(X,M)}\\
 \EA(M,X)&\EA(M)
\end{pmatrix}$ and $\begin{pmatrix}
\EA(M)&\widehat{\EA(M,Y)}\\\EA(Y,M)&\widehat{\EA(Y)}
\end{pmatrix}$ are subrings of $\EA(X\oplus M)$ and $\EA(M\oplus Y)$, respectively.


The following result is a general categorical property of the pair
$(A, B)$ where $A$ is a subring of $B$ with the same identity.
Each $B$-module can be regarded as an $A$-module just by the
restriction of scalars.

\begin{Lem}\label{Lem2}
Suppose that $A$ is a subring of $B$ with the same identity.

$(1)$ The restriction functor $F: \modcat{B} \ra \modcat{A}$ is an
exact faithful functor, and has a right adjoint $G=\Hom_A(_AB_B,-):
\modcat{A}\ra \modcat{B}$ and a left adjoint $E:=B\otimes_A-:
\modcat{A}\ra \modcat{B}$. In particular, $E$ preserves projective
modules and $G$ preserves injective modules.

$(2)$ The functor $E:=B\otimes_A-: \pmodcat{A}\ra \pmodcat{B}$ which
sends $A$ to $B$ is faithful.

$(3)$ For any $A$-module $M$ there is an $A$-homomorphism $\alpha_M:
GM\ra M$ such that the induced map $\Hom_B(X,GM)\ra\Hom_A(X,M)$ is
an isomorphism for all $B$-module $X$.
\end{Lem}

{\bf Proof.} We prove the statement $(2)$. The others are taken from
\cite[Lemma 4.2]{X1}.

Let $f: P\ra Q$ be a morphism in $\pmodcat{A}$. Since $Q$ is
projective as a left $A$-module, the morphism $\iota\otimes_A 1:
A\otimes_AQ\ra B\otimes_AQ$ is injective in $\pmodcat{A}$, where
$\iota$ is the inclusion map. Thus, it follows that $f=0$ from
$1\otimes f(p)=0$ with $p$ being an element in $P$.   $\square$
\medskip

Note that $\Lambda$ is a subring of $\Gamma$ with the same
identity. Then, by Lemma $\ref{Lem2}$, the functor
$\Gamma\otimes_\Lambda-: \Modcat{\Lambda}\ra\Modcat{\Gamma}$
preserves projective modules and is faithful, restricted to
$\pmodcat{\Lambda}$. Suppose that $\add^\Phi (V)$ is a full
subcategory of the $\Phi$-orbit category $D^b(\mathcal {A})^\Phi$
whose objects are the objects of $\add (V)$. It follows that the
additive functor $\Hom_{D^b(\mathcal {A})^\Phi}(V,-):
\add^\Phi(V)\ra \pmodcat{\Gamma}$ is an equivalence.
$$\xymatrix{
  \pmodcat{\Lambda}\ar[dr]_{G} \ar[r]^{\Gamma\otimes_\Lambda-} & \pmodcat{\Gamma}\ar[d]_F   \\
    & \add^\Phi(V)\ar@<-1.2ex>[u]_{\Hom_{D^b(\mathcal {A})^\Phi}(V,-)}  }$$

Let $F$ be the inverse of $\Hom_{D^b(\mathcal {A})^\Phi}(V,-)$, and
denote the composition of the functors $F$ and
$\Gamma\otimes_\Lambda-$ by $G:=F\circ(\Gamma\otimes_\Lambda-)$.
Thus the functor $G$ is faithful. We denote the image of $G$ by
$\mathcal {S}$ which is a subcategory of $\add^\Phi(V)$. The object
of $\mathcal{S}$ is $G(P)$ where $P$ is an object in
$\pmodcat{\Lambda}$. And the morphism set is defined as follows

$$\Hom_{\mathcal{S}}(G(P_1),G(P_2)):=(F\circ(\Gamma\otimes_{\Lambda}-))(P_1,P_2)$$
where $P_1, P_2$ are objects in $\pmodcat{\Lambda}$. Note that the
functor $\Gamma\otimes_\Lambda-$ is not full. Then the category
$\mathcal {S}$ is not a full subcategory of $\add^\Phi(V)$.

It is clear that the two categories $\pmodcat{\Lambda}$ and
$\mathcal {S}$ are equivalent under the additive functor $G$. Hence,
the category $\mathcal {S}$ is an additive subcategory of
$\add^\Phi(V)$. Note that $\Gamma\otimes_{\Lambda}\Lambda
e\cong\Gamma e$ as left $\Gamma$-modules where $e$ is an idempotent
of $\Lambda$. By Lemma $\ref{Lem2}$, we have that the functor
$\Gamma\otimes_{\Lambda}-$ is faithful. Thus, we have the following
sequence of isomorphisms, where short arguments are added on the
right hand side.
$$\begin{array}{lll}
\widehat{\EA(X,M)}&\cong&\Hom_{\Lambda}(\Lambda e_1,\Lambda e_2)\\
&\cong&\widehat{\Hom_{\Gamma}(\Gamma e_1, \Gamma e_2)}\quad(\mbox{
the functor }\Gamma\otimes_{\Lambda}- \mbox{ is faithful and
}\Gamma\otimes_{\Lambda}\Lambda e_i\cong \Gamma e_i)\\
&\cong&\Hom_{\mathcal {S}}(F(\Gamma e_1),F(\Gamma e_2))\quad
(F\mbox{ is a quasi-inverse of }\Hom_{D^b(\mathcal
{A})^{\Phi}}(V,-))\\
&\cong&\Hom_\mathcal{S}(X,M)
\end{array}$$
where $\widehat{\Hom_{\Gamma}(\Gamma e_1, \Gamma e_2)}$ is a subset
of $\Hom_{\Gamma}(\Gamma e_1, \Gamma e_2)$. Similarly, there exists
a set of $\mathbb{Z}$-module isomorphisms:
$$\begin{matrix}
 \Hom_\mathcal{S}(X,M)\cong\widehat{\EA(X,M)}&\Hom_\mathcal
{S}(M,Y)\cong\widehat{\EA(M,Y)}& \End_\mathcal
{S}(X)\cong\widehat{\EA(X)}\\
\End_\mathcal{S}(Y)\cong\widehat{\EA(Y)}& \Hom_\mathcal
{S}(X,Y)\cong\widehat{\EA(X,Y)}&\Hom_\mathcal
{S}(Y,X)\cong\EA(Y,X)\\
\End_\mathcal {S}(M)\cong\EA(M)&\Hom_\mathcal {S}(Y,M)\cong\EA(Y,M)&
\Hom_\mathcal {S}(M,X)\cong\EA(M,X).&
\end{matrix}$$

\medskip

Now, we prove the main result of this paper.

{\bf Proof of Theorem \ref{theorem}:} For simplicity of notation, we
will write $W, \overline{W}$ instead of $X\oplus M$ and $Y\oplus M$,
respectively. Define
$$\Lambda_1=\End_\mathcal {S}(W)=\begin{pmatrix}
  \widehat{\EA(X)}&\widehat{\EA(X,M)}\\
  \EA(M,X)&\EA(M)
\end{pmatrix},\quad\Lambda_2=\End_\mathcal{S}(\overline{W})=\begin{pmatrix}
\widehat{\EA(Y)}&\EA(Y,M)\\\widehat{\EA(M,Y)}&\EA(M)
\end{pmatrix}.$$

Since $0\ra X\stackrel{\alpha}\ra M_1\stackrel{\beta}\ra Y\ra 0$ is
an exact sequence in $\mathcal {A}$, we have an exact sequence
$$0\ra X\stackrel{\overline{\alpha}}\ra M_1\oplus M\stackrel{\overline{\beta}}\ra \overline{W}\ra 0$$
in $\mathcal {A}$, where $\overline{\alpha}:=(\alpha,0): X\ra
M_1\oplus M$ and
$\overline{\beta}:=\begin{pmatrix} \beta&0\\
0&1
\end{pmatrix}: M_1\oplus M\ra Y\oplus M$.

In order to prove Theorem \ref{theorem}, we first find out a tilting
complex over $\Lambda_1$.

\begin{Claim}\label{tiltingcomplex}
The complex
$$T^\bullet: 0\ra \Hom_\mathcal {S}(W,X)\stackrel{\Hom_\mathcal {S}(W,\overline{\alpha})}\lra
\Hom_\mathcal {S}(W,M_1\oplus M)\ra 0$$ is a tilting complex over
$\Lambda_1$.
\end{Claim}

{\bf Proof.} Note that $X$ and $M_1\oplus M$ belong to $\add(W)$.
Then both $\Hom_\mathcal {S}(W,X)$ and $\Hom_\mathcal
{C}(W,M_1\oplus M)$ are finitely generated projective left
$\Lambda_1$-modules. Thus, $T^\bullet$ is a complex in
$K^b(\pmodcat{\Lambda_1})$. It is immediate that the complex
$T^\bullet$ generates $K^b(\pmodcat{\Lambda_1})$ as a triangulated
category. It suffices to prove that $T^\bullet$ is self-orthogonal,
i.e.,
$\Hom_{K^b(\pmodcat{\Lambda_1})}(T^\bullet,T^\bullet[1])=\Hom_{K^b(\pmodcat{\Lambda_1})}(T^\bullet,T^\bullet[-1])=0.$

First, we will prove that
$\Hom_{K^b(\pmodcat{\Lambda_1})}(T^\bullet,T^\bullet[1])=0$. Note
that $W$ is an object in $\mathcal {S}$. Thus, there is an
equivalence of categories $\Hom_\mathcal {S}(W,-): \add^{\Phi}(W)\ra
\pmodcat{\Lambda_1}$. Let $s^\bullet=(s_i)_{i\in\Phi}$ be an element
in $\Hom_{K^b(\pmodcat{\Lambda_1})}(T^\bullet, T^\bullet[1])$. Then
there is a morphism $t=(t_i)_{i\in \Phi}\in \Hom_\mathcal
{S}(X,M_1\oplus M)=\widehat{\EA(X,M_1\oplus M)}$ satisfying that
$s^0=\Hom_\mathcal {S}(W,t)$. For simplicity, we denote
$\Hom_{\mathcal{S}}(X, Y)$ by $_{\mathcal{S}}(X, Y)$ in commutative
diagrams.
$$\xymatrix{
  & 0 \ar[d] \ar[r] & _\mathcal {S}(W,X) \ar[d]^{s^0=_\mathcal {S}(W,t)} \ar[r]^(.4){_\mathcal {S}(W,\overline{\alpha})}
  & _\mathcal {S}(W,M_1\oplus M)\ar@{-->}[ld]^{\widetilde{t}}\ar[d] \ar[r] & 0  \\
  0 \ar[r]& _\mathcal {S}(W,X) \ar[r]^{_\mathcal {S}(W,\overline{\alpha})} & _\mathcal {S}(W,M_1\oplus M) \ar[r] & 0 &   }\quad\quad\quad(\spadesuit)$$

By the definition of $\widehat{\EA(X,M)}$, there exists morphism
$\widetilde{t_i}: M_1\oplus M\ra (M_1\oplus M)[i]$ satisfying that
$t_i=\overline{\alpha}\widetilde{t_i}$ for $i\in \Phi$. For
abbreviate, we denote the morphism $(\widetilde{t_i})_{i\in \Phi}$
by $\widetilde{t}$. Now, we check that $\Hom_\mathcal
{S}(W,\overline{\alpha}) \widetilde{t}=\Hom_\mathcal {S}(W,t)$. Let
$(x_i)_{i\in \Phi}\in \Hom_\mathcal {S}(W,X)$. From the commutative
diagram $(\spadesuit)$, we have
$$\begin{array}{lll}
[\Hom_\mathcal {S}(W,\overline{\alpha}) \widetilde{t}]((x_i)_{i\in
\Phi})&=&\widetilde{t}((x_i\overline{\alpha}[i])_{i\in
\Phi})\\
&=&(\sum_{i,j\in \Phi \atop i+j=k}
x_i\overline{\alpha}[i]\widetilde{t_j}[i])_{k\in
\Phi}\\
&=&(\sum_{i,j\in \Phi \atop i+j=k} x_it_j[i])_{i\in
\Phi}\\
&=&\Hom_\mathcal {S}(W,t)((x_i)_{i\in \Phi}).\end{array}$$

Thus, the morphism $s^\bullet$ is null-homotopic, that is,
$\Hom_{K^b(\pmodcat{\Lambda_1})}(T^\bullet,T^\bullet[1])=0$.

Next, we will prove that
$\Hom_{K^b(\pmodcat{\Lambda_1})}(T^\bullet,T^\bullet[-1])=0$. Let
$s^\bullet$ be an element in
$\Hom_{K^b(\pmodcat{\Lambda_1})}(T^\bullet,T^\bullet[-1])$. Then
there is a morphism $t'$ in $\Hom_{\mathcal{S}}(M_1\oplus M, X)$
such that $s^0=\Hom_{\mathcal{S}}(W, t')$.
$$\xymatrix{
  0 \ar[r]& _\mathcal {S}(W,X) \ar[d] \ar[r]^(.4){_\mathcal {S}(W,\overline{\alpha})} & _\mathcal {S}(W,M_1\oplus M)
  \ar[d]^{s^0=_\mathcal {S}(W,t')} \ar[r]& 0 \ar[d] & \\
  & 0 \ar[r] & _\mathcal {S}(W,X) \ar[r]^{_\mathcal {S}(W,\overline{\alpha})} & _\mathcal {S}(W, M_1\oplus M) \ar[r] & 0   }\quad\quad(\spadesuit1)$$
By the commutativity of $(\spadesuit1)$, we have $\Hom_\mathcal
{S}(W,t')\Hom_\mathcal {S}(W,\overline{\alpha})=\Hom_\mathcal
{S}(W,t'\overline{\alpha})=0$. This implies that
$t\overline{\alpha}=0$ in $\End_\mathcal {S}(M_1\oplus M)$. So we
have $t_i\overline{\alpha}[i]=0$ for $i\in\Phi$. Note that the
morphism $\overline{\alpha}$ is monic in $\mathcal {A}$, we can
deduce $t_0=0$ in $\mathcal {A}$. Since $\Ext^i_\mathcal {A}(M,X)=0$
for $0\neq i\in\Phi$, we get $t'_i=0$ for $0\neq i\in\Phi$. Hence,
$\Hom_{K^b(\pmodcat{\Lambda_1})}(T^\bullet,T^\bullet[-1])=0$.$\square$
\medskip

Now, we will determine the endomorphism ring of $T^\bullet$.
\begin{Claim}\label{endomorphism} The two
rings $\End_{K^b(\pmodcat{\Lambda_1})}(T^\bullet)$ and $\Lambda_2$
are isomorphic.
\end{Claim}

{\bf Proof.} Let $f^\bullet$ be a morphism in
$\End_{K^b(\pmodcat{\Lambda_1})}(T^\bullet)$. We have the following
commutative diagram
$$\xymatrix{
  0  \ar[r] & _\mathcal {S}(W,X) \ar[d]_{f^0=_{\mathcal{S}}(W, u)} \ar[r]^{_\mathcal {S}(W,\overline{\alpha})} & _\mathcal {S}(W,M_1\oplus M)
  \ar[d]_{f^1=_{\mathcal{S}}(W, v)} \ar[r] & 0  \\
  0 \ar[r] & _\mathcal {S}(W,X) \ar[r]^{_\mathcal {S}(W,\overline{\alpha})} & _\mathcal {S}(W,M_1\oplus M) \ar[r] & 0 .  }\quad\quad(\spadesuit2)$$

Note that there is an equivalence of categories $\Hom_\mathcal
{S}(W,-): \add (W)\ra \pmodcat{\Lambda_1}$. Then there exist
morphisms $u=(u_i)_{i\in\Phi}\in \End_\mathcal {S}(X)\cong
\widehat{\EA(X)}, v=(v_i)_{i\in\Phi}\in\End_\mathcal {S}(M_1\oplus
M)\cong \EA(M_1\oplus M)$ such that $f^0=\Hom_\mathcal {S}(W,u),
f^1=\Hom_\mathcal {S}(W,v)$. By the commutative diagram
$(\spadesuit2)$, we have $\Hom_\mathcal
{S}(W,\overline{\alpha}v)=\Hom_\mathcal {S}(W,u\overline{\alpha})$.
It implies that $\overline{\alpha}v=u\overline{\alpha}$, i.e.,
$\overline{\alpha}v_i=u_i\overline{\alpha}[i]$ for $i\in\Phi$. For
each pair of $(u_i, v_i)$, there exists a commutative diagram
$$\xymatrix{
  X \ar[d]_{u_i} \ar[r]^{\overline{\alpha}} & M_1\oplus M \ar[d]_{v_i} \ar[r]^{\overline{\beta}} & Y\oplus M \ar@{-->}[d]_{h_i} \ar[r]^{\overline{w}}
  & X[1] \ar[d]^{u_i[1]} \\
  X[i] \ar[r]^{\overline{\alpha}[i]} & (M_1\oplus M)[i] \ar[r]^{\overline{\beta}[i]} & (Y\oplus M)[i] \ar[r]^{\overline{w}[i]} & X[i+1]   }\quad\quad(\spadesuit3)$$
for $i\in \Phi$. Thus there is a morphism $h_i:{}Y\oplus
M\ra(Y\oplus M)[i]$ such that
$\overline{\beta}h_i=v_i\overline{\beta}[i]$ for $i\in\Phi$. We
denote $(h_i)_{i\in\Phi}$ by $h$, and denote $h_i, v_i$ by
$\begin{pmatrix}h_{i1}&h_{i2}\\h_{i3}&h_{i4}
\end{pmatrix}, \begin{pmatrix}
  v_{i1}&v_{i2}\\v_{i3}&v_{i4}
\end{pmatrix}$, respectively. It follows that $h_{i1}=v_{i1}\beta[i],
h_{i3}=v_{i3}\beta$ from
$\overline{\beta}h_i=v_i\overline{\beta}[i]$ for $i\in \Phi$. By the
definition of $\widehat{\EA(Y\oplus M)}$, we deduce that $h$ is an
element in $\End_\mathcal {S}(Y\oplus M)$.

Thus, we can define the following correspondence:
$$\Psi: \End_{K^b(\pmodcat{\Lambda})}(T^\bullet)\ra \End_\mathcal {S}(Y\oplus M)$$
$$f^\bullet\longmapsto h=(h_i)_{i\in \Phi}.$$

First, we will prove that the map $\Psi$ is well-defined. Suppose
that $f^\bullet$ is null-homotopic. Thus there is a morphism $s:
\Hom_\mathcal {S}(W,M_1\oplus M)\ra \Hom_\mathcal {S}(W, X)$ such
that $f^0=\Hom_\mathcal {S}(W,\overline{\alpha})s,
f^1=s\Hom_\mathcal {S}(W,\overline{\alpha})$. We denote $s$ by
$\Hom_\mathcal {S}(W,t)$ with $t=(t_i)_{i\in \Phi}\in \Hom_\mathcal
{S}(M_1\oplus M, X)$. By \cite[Lemma 3.5]{HX2}, we can get
$u_i=\overline{\alpha}t_i, v_i=t_i\overline{\alpha}[i]$ for $i\in
\Phi$. Thus, by the commutative diagram $(\spadesuit3)$, we have
$$\overline{\beta}h_i=v_i\overline{\beta}[i]=t_i\overline{\alpha\beta}[i]=0, h_i\overline{w}[i]=\overline{w}u_i[1]=0$$
for $i\in\Phi$. So there are morphisms $l_i: Y\oplus M\ra (M_1\oplus
M)[i], p_i: X[1]\ra (M_1\oplus M)[i]$ such that
$h_i=l_i\overline{\beta}[i], h_i=\overline{w}p_i$ for $i\in\Phi$.
Since $\overline{w}=0$, restricted to $M$, and $\Ext^i_A(Y,M)=0$ for
$i\in \Phi$, we have $h_i=0$ for $0\neq i\in \Phi$. Note that
$\mathcal {A}$ can be embedded into $D^b(\mathcal {A})$. Since
$\overline{\beta}, h_0$ are morphisms in $\mathcal {A}$ and
$\overline{\beta}$ is an epimorphism in $\mathcal {A}$, we can
deduce $h_0=0$ in $\mathcal {A}$. Hence the map $\Psi$ is
well-defined.

Second, we will prove that the map $\Psi$ is injective. Suppose
$h=0$ in $\mathcal {S}$. Then we have $v_i\overline{\beta}[i]=0,
\overline{w}u_i[1]=0$ for $i\in \Phi$. So, by the properties of
cohomological functor, we have that $u_i, v_i$ factorize through
$\overline{\alpha}$ and $\overline{\alpha}[i]$, respectively, i.e.,
there are morphisms $a_i: M_1\oplus M\ra X[i], b_i: M_1\oplus M\ra
X[i]$ such that $v_i=a_i\overline{\alpha}, u_i=\overline{\alpha}b_i$
for $i\in \Phi$. It follows that $u_i=0=v_i$ from $\Ext_\mathcal
{A}^i(M,X)=0$ for $0\neq i\in\Phi$. When $i=0$, we can get
$u_0\overline{\alpha}=\overline{\alpha}v_0=\overline{\alpha}b_0\overline{\alpha}$.
Note that $u_0, b_0$ are morphisms in $\mathcal {A}$ and
$\overline{\alpha}$ is injective in $\mathcal {A}$, it follows that
$u_0=\overline{\alpha}b_0$. Set $\iota_0(a_0): \Hom_\mathcal
{S}(M_1\oplus M, X)\ra \EA(M_1\oplus M, X)$, where $\iota_0(a_0)$ is
an element in $\EA(M_1\oplus M, X)$ concentrated in degree $0$. It
is easy to check $f^1=\mu(\iota_0(a_0)))\Hom_\mathcal
{S}(W,\overline{\alpha})$. Since
$\overline{\alpha}a_0\overline{\alpha}=\overline{\alpha}v_0=u_0\overline{\alpha}$,
and $\overline{\alpha}$ is monic in $\mathcal {A}$, we have
$\overline{\alpha}a_0=u_0$. Thus, we have $f^0=\Hom_\mathcal
{C}(W,\overline{\alpha})\mu(\iota_0(a_0)))$. Altogether, we have
proved that the morphism $f^\bullet$ is null-homotopic. Hence, the
map $\Psi$ is injective.

Third, we will prove that the map $\Psi$ is surjective. Let $h$ be a
morphism in $\Lambda_2$. Then $h$ can be written as the form
$h=(h_i)_{i\in\Phi}$ with $h_i=\begin{pmatrix}
h_{i1}&h_{i2}\\h_{i3}&h_{i4}
\end{pmatrix}: Y\oplus M\ra (Y\oplus M)[i]$ in $D^b(\mathcal {A})$. By the
definition of $\widehat{\EA(M,Y)}$, there exists a morphism
$\begin{pmatrix}s_{i1}\\
s_{i2}\end{pmatrix}: M_1\oplus M\ra M_1[i]$ such that
$\begin{pmatrix}s_{i1}\\
s_{i2}\end{pmatrix}\beta[i]=\overline{\beta}\begin{pmatrix}
h_{i1}\\h_{i3}
\end{pmatrix}$ for $i\in\Phi$. Thus there exists a morphism
$\begin{pmatrix}
  s_{i1}&\beta h_{i2}[-i]\\s_{i2}&h_{i4}[-i]
\end{pmatrix}: M_1\oplus M\ra (M_1\oplus M)[i]$ satisfying that $\begin{pmatrix}
  s_{i1}&\beta
  h_{i2}[-i]\\s_{i2}&h_{i4}[-i]
\end{pmatrix}\overline{\beta}[i]=\overline{\beta}h_i$ for
$i\in\Phi$. We denote $\begin{pmatrix}s_{i1}&\beta
h_{i2}[-i]\\s_{i2}&h_{i4}[-i]
\end{pmatrix}$ by $v_i$ for $i\in \Phi$. So there exist $u_i: X\ra X[i]$ such that
$\overline{\alpha}v_i=u_i\overline{\alpha}[i]$ for $i\in \Phi$. Let
us denote $(u_i)_{i\in\Phi}$ and $(v_i)_{i\in\Phi}$ by $u$ and $v$,
respectively. Since
$(\overline{\alpha}v)_{i\in\Phi}=\overline{\alpha}v_i=u_i\overline{\alpha}[i]=(u\overline{\alpha})_i$
for $i\in \Phi$, we have $\Hom_\mathcal
{C}(W,\overline{\alpha})\Hom_\mathcal {S}(W,v)=\Hom_\mathcal
{C}(W,u)\Hom_\mathcal {S}(W,\overline{\alpha})$. So the map $\Psi$
is surjective.

At last, we will prove that the map $\Psi$ is a ring homomorphism.
Let $f^\bullet, g^\bullet$ be morphisms in $\End_\mathcal
{C}(T^\bullet)$. Suppose that $f^0=\Hom_\mathcal
{C}(W,u),f^1=\Hom_\mathcal {S}(W,v)$, $g^0=\Hom_\mathcal {S}(W,u'),
g^1=\Hom_\mathcal {S}(W,v')$ with $u,u'\in \E_\mathcal {S}(X),
v,v'\in \E_\mathcal {S}(M_1\oplus M)$. Set $\Psi(f^\bullet)=h$,
$\Psi(g^\bullet)=h'$. By computation, we can get
$uu'=(\sum_{i+j=k\atop i, j\in \Phi}u_iu_j'[i])_{k\in \Phi},
vv'=(\sum_{i+j=k\atop  i, j\in \Phi}v_iv_j'[i])_{k\in \Phi}$, and
$\Psi(f^\bullet g^\bullet)=(\sum_{i+j=k\atop i,j\in
\Phi}h_ih_j'[i])_{k\in \Phi}$. For $k\in \Phi$,
$$\overline{\beta}(\sum_{i+j=k\atop i,j\in
\Phi}h_ih_j'[i])=\sum_{i+j=k\atop i,j\in
\Phi}\overline{\beta}h_ih_j'[i]=\sum_{i+j=k\atop i,j\in
\Phi}v_i\overline{\beta}[i]h_j'[i]=\sum_{i+j=k\atop i,
j\in\Phi}v_iv_j'[i]\overline{\beta}[k],$$
$$\overline{w}\{(\sum_{i+j=k\atop i,j\in\Phi}u_iu_j'[i])[1]\}=\sum_{i+j=k\atop i,j\in\Phi}\overline{w}u_i[1]u_j'[i+1]
=\sum_{i+j=k\atop
i,j\in\Phi}h_i\overline{w}[i]u_j'[i+1]=\sum_{i+j=k\atop
i,j\in\Phi}h_ih_j'[i]\overline{w}[k].$$

It follows that $\Psi(f^\bullet
g^\bullet)=\Psi(f^\bullet)\Psi(g^\bullet)$. Hence the map $\Psi$ is
a ring homomorphism. $\square$

\medskip
If $\Phi=\{0\}$, we have the following corollary which suggests that
any exact sequence implies a derived equivalence between two rings.

\begin{Koro}\label{excat}Let $\mathcal {A}$ be an abelian category and $M$ an object in
$\mathcal{A}$. Suppose that
$$0\ra X\stackrel{\alpha}\ra M_1\stackrel{\beta}\ra Y\ra 0$$
is an exact sequence in $\mathcal{A}$ with $M_1\in\add(M)$.

Then the two rings
$$\begin{pmatrix}
  \widehat{\End_{\A}(X)}&\widehat{\Hom_{\A}(X,M)}\\
  \Hom_{\A}(M,X)&\End_{\A}(M)
\end{pmatrix}\mbox{ and }  \begin{pmatrix}\widehat{\End_{\A}(Y)}
  &\Hom_{\A}(Y,M)\\\widehat{\Hom_{\A}(M,Y)}
  &\End_{\A}(M)
\end{pmatrix}$$
are derived equivalent, where
$$\widehat{\End_{\A}(Y)}:=\{t\in\End_{\A}(Y)\mid \mbox{there exists a morphism }t_1:{}M_1\ra M_1 \mbox{ such that }
\beta t=t_1\beta\mbox{ in } \A\}$$
$$\widehat{\End_{\A}(X)}:=\{t\in\End_{\A}(X)\mid \mbox{there exists a morphism }t_1:{}M_1\ra M_1 \mbox{ such that }
t\alpha=\alpha t_1 \mbox{ in } \A\}$$
$$\widehat{\Hom_{\A}(X,M)}:=\{t\in\Hom_{\A}(X,M)\mid \mbox{there exists a morphism }t_1:{}M_1\ra M \mbox{ such that }
t=\alpha t_1\mbox{ in } \A\}$$
$$\widehat{\Hom_{\A}(M,Y)}:=\{t\in\Hom_{\A}(M,Y)\mid \mbox{there exists a morphism }t_1:{}M\ra M_1 \mbox{ such that }
t=t_1\beta\mbox{ in } \A\}.$$
\end{Koro}

\medskip

\begin{Koro}\label{cor3}Let $A$ and $B$ be two rings with the same identity and $\alpha:{}B\hookrightarrow A$ be an injective ring
homomorphism. Then the two rings
$$\End_B(B\oplus{}_BA) \mbox{ and }
\begin{pmatrix}\widehat{\End_B(A/B)}&\Hom_B(A/B,{}_BA)\\\widehat{\Hom_B(_BA,A/B)}
&\End_B(A)
\end{pmatrix}$$ are derived equivalent. In particular, if $A$ is finitely generated projective as a left $B$-module. Then the
two rings $B$ and
$\begin{pmatrix}\End_B(A/B)&\Hom_B(A/B,A)\\\Hom_B(A,A/B) &\End_B(A)
\end{pmatrix}$
are derived equivalent.
\end{Koro}

{\bf Proof.} Any left $B$-homomorphism $f: B\ra{} _BA$ is determined
by $f(1)$ where $1$ is the identity of $B$. Let $g:{}A\ra A$ be the
a left $A$-module homomorphism which sends $1$ to $f(1)$. Then $g$
is a left $B$-module and satisfies $f=\alpha g$. Thus, $\alpha$ is a
left $\add({}_BA)$-approximation of B in $\Modcat{B}$. If $_BA$ is
finitely generated projective, then $B\oplus{}_BA$ is a projective
generator over $B$. So the two rings $B$ and $\End_B(B\oplus{}_BA)$
are Morita equivalent. Clearly, the morphism $\beta$ is a right
$\add({}_BA)$-approximation of $A/B$. Thus, the two rings $B$ and
$\begin{pmatrix}\End_B(A/B)&\Hom_B(A/B,A)\\\Hom_B(A,A/B) &\End_B(A)
\end{pmatrix}$ are derived equivalent.
$\square$

\medskip

Recall that a triangulated category $\mathcal {T}$ is said to be
{\em algebraic} if it is triangle equivalent to the stable category
of an exact Frobenius category $(\mathcal {B}, S)$ where $\mathcal
{B}$ denotes an extension closed full subcategory of an abelian
category $\mathcal {A}$ and $S$ is the set of exact sequences in
$\mathcal {A}$ with terms in $\mathcal {B}$.

\begin{Koro}\label{D} Let $\Phi$ be an admissible subset of
$\mathbb{N}$. Let $\mathcal {T}$ be an algebraic triangulated
category and $M$ be an object in $\mathcal{T}$. Suppose that
$X\stackrel{\alpha}\ra M_1\stackrel{\beta}\ra Y\ra \Sigma X$ is a
triangle in $\mathcal {T}$ where $M_1\in \add (M)$. Then the two
rings
$$\begin{pmatrix}
  \widehat{{\rm E}^\Phi_\mathcal {T}(X)}&\widehat{{\rm E}^\Phi_\mathcal {T}(X,M)}\\
  {\rm E}^\Phi_\mathcal {T}(M,X)&{\rm E}^\Phi_\mathcal {T}(M)
\end{pmatrix}\mbox{ and }
\begin{pmatrix} \widehat{{\rm E}^\Phi_\mathcal {T}(Y)}&{\rm E}^\Phi_\mathcal {T}(Y,M)\\\widehat{{\rm E}^\Phi_\mathcal {T}(M,Y)}&{\rm E}^\Phi_\mathcal {T}(M)
\end{pmatrix}$$
are derived equivalent.
\end{Koro}

{\bf Proof.} The proof is similar to Theorem $\ref{theorem}$. By
\cite[Lemma 3.1 ]{D}, the map $\Psi$ is well-defined and injective
without any condition when $\mathcal {T}$ is an algebraic
triangulated category.  $\square$

\medskip

Let $\Phi=\{0\}$ or $\mathbb{Z}$. If $\alpha$ is a left
$\add(M)$-approximation of $X$ in $\mathcal {T}^{\Phi}$ and if
$\beta$ is a right $\add(M)$-approximation of $Y$ in $\mathcal
{T}^{\Phi}$, then we can cover the main result of \cite{D}.

\medskip

%
%
%

Let $A$ be a ring. We denote by $\gld(A)$ the left global dimension
of $A$. For a left $A$-module $M$, $\injdim({}_AM)$ (respectively,
$\pd(_{}AM)$) means the left injective (respectively, projective)
dimension of $_AM$. The big finitistic dimension, $\Fd(A)$, is
defined as follows.
$$\Fd(A)={\rm sup}\{\pd({}_AM)\mid M\in \Modcat{A} \mbox{ and }\pd(M)<\infty\}.$$
Similarly, the small finitistic dimension, $\fd(A)$, is defined as
$$\fd(A)={\rm sup}\{\pd({}_AM)\mid M\in \modcat{A} \mbox{ and }\pd(M)<\infty\}.$$

In \cite{K}, Kato discussed the difference between the left global
dimension of two rings which are derived equivalent. In \cite{PX},
Pan and Xi proved that the difference between the small finitistic
dimension of two coherent rings is less than the term length of a
tilting complex. Combining with Theorem \ref{theorem}, we have the
following Corollary.

\begin{Koro} Keep the notation as in Theorem $\ref{theorem}$. Set
$$\Lambda_1:=\begin{pmatrix}
  \widehat{\EA(X)}&\widehat{\EA(X,M)}\\
  \EA(M,X)&\EA(M)
\end{pmatrix}\mbox{ and } \Lambda_2:=\begin{pmatrix} \widehat{{\rm E}^\Phi_A(Y)}&{\rm E}^\Phi_A(Y,M)\\\widehat{{\rm E}^\Phi_A(M,Y)}&{\rm E}^{\Phi}_A(M).
\end{pmatrix}$$
Then the following statements hold.

$(1)$ $\gld(\Lambda_1)-1\leq\gld(\Lambda_2)\leq\gld(\Lambda_1)+1$.

$(2)$
$\injdim(_{\Lambda_1}\Lambda_1)-1\leq\injdim(_{\Lambda_2}\Lambda_2)\leq\injdim(_{\Lambda_1}\Lambda_1)+1$.

$(3)$ Suppose that $\Lambda_1$ and $\Lambda_2$ are left coherent
rings. Then
$$\fd(\Lambda_1)-1\leq\fd(\Lambda_2)\leq\fd(\Lambda_1)+1.$$
\end{Koro}

{\bf Proof.} By Theorem $\ref{theorem}$, \cite[Proposition 1.7]{K}
and \cite[Theorem 1.1]{PX}, we can get the conclusion. $\square$
\medskip

\medskip

\section{Example}

In this section, we will give an example to illustrate our results.

Let $A$ be the following $k$-algebra with quiver

$$
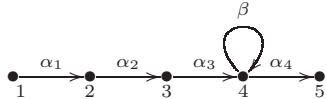
\xymatrix{*{\bullet}\ar[r]_(0.1){1}_(1){2}^{\alpha_1}&*{\bullet}\ar[r]_(0.1){}_(1){3}^{\alpha_2}&*{\bullet}\ar[r]_(0.1){}_(1){4}^{\alpha_3}
&*{\bullet}\ar[r]_(0.1){}_(1){5}^{\alpha_4}\ar@(ul,ur)^{\beta}
&*{\bullet} }$$ satisfying that $\rad^2(A)=0$.
\medskip

%

Let $M$ be \unitlength=1.00mm \special{em:linewidth 0.4pt}
\linethickness{0.4pt}
\begin{picture}(8.00,7.00)
\put(0,3){$3$} \put(2,-1){$4$}  \put(4,3){$4$}
\end{picture} which is the injective envelop of $S_4$. Thus, we have the following sequence
$$0\ra \Omega(M)\ra P(M)\ra M\ra 0$$ where $P(M)=$ \unitlength=1.00mm \special{em:linewidth 0.4pt}
\linethickness{0.4pt}
\begin{picture}(15.00,7.00)
\put(0,-1){$4$}\put(0,3){$3$}\put(4,1){$\oplus$} \put(8,-1){$4$}
\put(10,3){$4$} \put(12,-1){$5$}
\end{picture}
is the projective cover of $M$, and $\Omega(M)=4\oplus 5$ is the
first syzygy of $M$. Then $\End_A(\Omega(M)\oplus P(M))$ and
$\widehat{\End_A(\Omega(M)\oplus P(M))}$ can be described by quiver
with relations:
%
%
%
%
%
%
%

$$\begin{array}{ll}
\End_A(\Omega(M)\oplus
P(M))\qquad\qquad\qquad&\widehat{\End_A(\Omega(M)\oplus
P(M))}\\
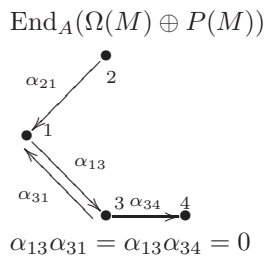
\xymatrix{&*{\bullet}\ar[ld]^(0.1){2}^(1){}_{\alpha_{21}}&\\
*{\bullet}\ar[dr]^(0.1){1}^(1){3}^{\alpha_{13}}&&\\
&*{\bullet}\ar@<4pt>[ul]^{\alpha_{31}}\ar[r]^(0.1){}^(1){4}^{\alpha_{34}}&*{\bullet}
}&
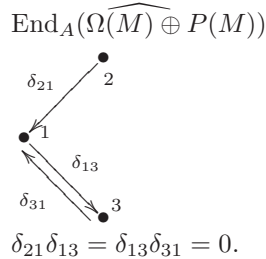
\xymatrix{&*{\bullet}\ar[ld]^(0.1){2}^(1){}_{\delta_{21}}\\
*{\bullet}\ar[dr]^(0.1){1}^(1){3}^{\delta_{13}}&\\
&*{\bullet}\ar@<4pt>[ul]^{\delta_{31}}}\\
\alpha_{13}\alpha_{31}=\alpha_{13}\alpha_{34}=0&\delta_{21}\delta_{13}=\delta_{13}\delta_{31}=0.\end{array}$$

The endomorphism algebra $\End_A(M\oplus P(M))$ is given by the
following quiver

$$
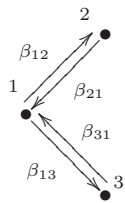
\xymatrix{&*{\bullet}\ar[dl]^{\beta_{21}}\\
*{\bullet}\ar@<4pt>[ur]^(0.1){1}^(1){2}^{\beta_{12}}\ar[dr]^(0.1){}^(1){3}_{\beta_{13}}&\\
&*{\bullet}\ar@<-4pt>[ul]^(0.1){}^(1){}_{\beta_{31}}}$$ with
relations $\{\beta_{12}\beta_{21}=\beta_{13}\beta_{31},
\beta_{21}\beta_{12}, \beta_{31}\beta_{12}\}$. Note that $P(M)$ is a
finitely generated projective left $A$-module. By Corollary
\ref{excat}, we have that the two rings
$\widehat{\End_A(\Omega(M)\oplus P(M))}$ and $\End_A(M\oplus P(M))$
are derived equivalent.

\medskip \footnotesize{
}

\end{document}